\documentclass[preprint,12pt]{elsarticle}
\usepackage{amsfonts}
\usepackage{amsmath}
\usepackage{mathrsfs}
\usepackage[applemac]{inputenc}
\usepackage{amssymb}
\usepackage{times, ulem, amsfonts}
\usepackage{mathpazo, amsmath, amssymb, amsthm,color}
\makeatletter
\def\ps@pprintTitle{%
   \let\@oddhead\@empty
   \let\@evenhead\@empty
   \let\@oddfoot\@empty
   \let\@evenfoot\@oddfoot
}
\makeatother
\def\y{{\mathcal Y}}
\def\z{{\mathcal X}}

\newcommand{\R}{{\mathbb{R}}}
\newcommand{\Z}{{\mathbb{Z}}}

\def\div{ \hbox{\rm div}\,  }
\def\ta{\vartheta}
\def\La{\Lambda}

\def\ddj{\dot \Delta_j}

\newtheorem{theorem}{Theorem}[section]
\newtheorem{lemma}[theorem]{Lemma}
\newtheorem{definition}[theorem]{Definition}

\newtheorem{proposition}[theorem]{Proposition}
\newtheorem{remark}[theorem]{Remark}
\numberwithin{equation}{section}

\begin{document}
\title{ 
 Long-time behavior for      three dimensional  compressible
viscous and heat-conductive gases}

  \author{Xiaoping Zhai}
\ead{zhaixp@szu.edu.cn}

 \author{Zhi-Min Chen\corref{cor1}}
 \ead{zmchen@szu.edu.cn}
 \address{School  of Mathematics and Statistics, Shenzhen University, Shenzhen 518060, China}
\cortext[cor1]{Corresponding author}

\baselineskip=24pt

\begin{abstract}
Large-time behavior of solutions to the compressible Navier-Stokes equations
for  viscous and heat-conductive gases in $\mathbb{R}^3$ is examined.
Under  a suitable  condition involving only the low frequencies of the  initial data, optimal time decay rates  for  the  non-isentropic compressible Navier-Stokes flows are obtained, by developing some  energy arguments given by  [Xin and Xu, ArXiv:1812.11714v2].
\end{abstract}

\maketitle


\setcounter{section}{0}
\setcounter{theorem}{0}
\section{Introduction and the main result}
In this paper, we consider the  three-dimensional
  non-isentropic compressible Navier-Stokes equations:
\begin{eqnarray}\label{mm}
\left\{\begin{aligned}
&\partial_t \rho + \div(\rho u) = 0\, ,\\
 &     \partial_t ( \rho u )+ \div ( \rho u \otimes u )-\mu\Delta u-(\lambda+\mu)\nabla\div u + \nabla P = 0,\\
      &\rho (\partial_t\theta+u\cdot\nabla\theta)
+ P\div u-\kappa\Delta\theta=2\mu |D(u)|^2+\lambda(\div u)^2,\\
    &(\rho,u,\theta)|_{t=0}=(\rho_0,u_0,\theta_0)
\end{aligned}\right.
\end{eqnarray}
for the velocity   $u = (u_1, u_2, u_3),$ the density  $\rho,$ the absolute temperature  $\theta$,  the pressure
$P$   and the deformation tensor  $D(u) = \frac{1}{2} \big( \nabla u + (\nabla u)^{T} \big)$.
The  viscosity coefficients $\mu$ and $\lambda$ are subject to  the standard  physical restrictions:
\begin{align*}
\mu>0\quad\hbox{and}\quad
3\lambda+2\mu>0.
\end{align*}
 $\kappa>0$ denotes  the heat
conductivity coefficient.
For a general fluid, $P$ is a function of  $\rho$ and $\theta$.
Here, we only consider the case of perfect heat conducting and viscous gases, i.e., $P(\rho,\theta)=A\rho\theta$ for some constant $A>0$.
 Without loss of generality, we assume that $A=1$ throughout the paper.

If the effect of the temperature is neglected and thus  the pressure is a  function of $\rho$,  Eq. (\ref{mm}) reduces to the isentropic  compressible
Navier-Stokes equations, which have been widely studied (see \cite{charve, chenqionglei, chenqionglei3, chenqionglei2, zhaixiaoping, danchin2000, helingbing, xujiang2017arma, Feireisl, Feireisl3, Feireisl2, haspot, huangxiangdi2, xujiang2019cmp} and the references therein).
However, the system (\ref{mm}), involving the heat conducting effect, is  much more complicated due to the occurrence of  stronger nonlinearity, and thus  has attracted increasing  attention from  mathematicians.
  Significant  progress has been made in  the understanding of   the existence,  regularity and asymptotic behavior of the solutions over the past several decades.

 In the absence of vacuum,  the initial density $\rho_0$ is bounded away from zero.  Nash \cite{nash} obtained the local existence of classical  solutions in Sobolev spaces.
Itaya  \cite{itaya}  showed the existence of  the local classical solutions in H\"older spaces, and this result was further derived by   Tani \cite{tani} in a domain with boundary. Later on, Matsumura and Nishida \cite{ma2} proved the existence of  global  small solutions and stability of a basic steady-state solution ( see also \cite{ks} for the one-dimensional situation). For the existence and asymptotic behavior of  weak solutions,  Jiang \cite{jiangsong1998} considered the one-dimensional case, Jiang \cite{jiangsong1999} examined  a spherically symmetric structure  in two- and three-dimensional spaces,  Hoff and Jenssen \cite{hoff} discussed  spherically and cylindrically symmetric structures  in the three-dimensional space, and Feireisl \cite{Feireisl3}  investigated the problem  in two- and three-dimensional bounded domains.

In the presence of  vacuum, long-time  solutions  may not exist. If $ \kappa= 0$, it is obtained  by  Xin \cite{xinzhouping} that any non-zero smooth
solution with initial compactly supported density would blow up in a finite time.
 Feireisl \cite{Feireisl} got  the existence of  variational solutions in $\mathbb{R}^d\ (d\ge 2)$ under the condition that the temperature equation  is satisfied only as an inequality in the sense of distributions. For the  equality equations  in the sense of distributions, Bresch and Desjardins \cite{bresch} proposed some different assumptions from \cite{Feireisl} and obtained the existence of global weak solutions to the  compressible Navier-Stokes equations with large initial data in $T^3$ and  $\mathbb{R}^3$. Huang and Li \cite{huangxiangdi4} established the global existence and uniqueness of classical solutions to the three-dimensional  compressible Navier-Stokes system with smooth initial data,
which are of small energy but possibly large oscillations, and with  initial density allowed to vanish. Wen and Zhu \cite{wenhuanyao} established the global existence of spherically and cylindrically symmetric classical  solutions to the three-dimensional  compressible Navier-Stokes equations.  We also emphasis some blowup criterions given by Huang and Li  \cite{huangxiangdi2, huangxiangdi4,huangxiangdi1} and Huang {\it et al.} \cite{huangxiangdi3}  (see also the references therein). The readers may refer to Chen {\it et al.} \cite{chenqionglei} and Wen and Zhu \cite{wenhuanyao} for more recent developments  on this subject.

Let us recall the  critical  regularity  framework of \eqref{mm}.
Rigorously speaking, Eq.  \eqref{mm} does not possess any scaling invariance with respect to $(u,\rho, \theta, P)$. However, in the absence of   the pressure,  we define
 the scaling transformation:
\begin{equation}\label{scaling}
\rho(t,x)\rightarrow\rho(\lambda^2t,\lambda x),\quad
u(t,x)\rightarrow \lambda u(\lambda^2t,\lambda x), \quad  \theta(t,x) \rightarrow\lambda^2\theta(\lambda^2t,\lambda x)  \qquad\lambda>0.
\end{equation}
A function space is said to be  critical with respect to (\ref{mm})  if  the space norm is invariant with respect to the scaling
transformation (\ref{scaling}).
For example,  the product space
$\dot{B}_{p,1}^{\frac {d}{p}}({\mathbb R} ^d)
\times \dot{B}_{p,1}^{\frac {d}{p}-1}({\mathbb R} ^d) \times \dot{B}_{p,1}^{\frac {d}{p}-2}({\mathbb R} ^d) $,
$1\le p\le \infty$,
is a critical space for the system (\ref{mm}).

In the  framework of  critical spaces, a breakthrough was made by
Danchin \cite{danchin2000} for the  isentropic  compressible
Navier-Stokes equations, where he proved the local wellposedness with large initial data and global solutions with small initial data.
This result was further extended by  Charve {\it et al.} \cite{charve},
 Chen {\it et al.} \cite{chenqionglei}, Haspot \cite{haspot} and the authors \cite{zhaixiaoping}. 
 For the non-isentropic compressible Navier-Stokes equations \eqref{mm}, Chikami and Danchin \cite{chikami} and Danchin  \cite{danchin2001cpde} considered the  local wellposedness problem.  Global small solutions in a critical  $L^p$ Besov  framework were  obtained respectively  by Danchin \cite{danchin2001arma} for  $p=2$ and Danchin and He \cite{helingbing} for more general $p$. It should be mentioned that the critical Besov space  used in \cite{chikami} and \cite{helingbing}
seems the largest one in which the system (\ref{mm}) is well-posed.
Indeed, Chen {\it et al.} \cite{chenqionglei2} proved the ill-posedness of (\ref{mm})
in $\dot B^{\frac  3p}_{p,1}\times \dot B^{\frac 3p-1}_{p,1}\times \dot B^{\frac 3p-2}_{p,1}(\R^3)$ for $p>3$.

In the  present paper,
we consider  the large time behaviour of  three-dimensional global solutions
  to  \eqref{mm}, stemmed from  the pioneering work of  Matsumura and Nishida \cite{ma2} on the decay estimate
\begin{eqnarray}\label{R-E7}
\sup_{t\geq0}\: \langle t\rangle^{\frac{3}{4}} \|(\rho, u,\theta)(t)-(\bar{\rho},0,\bar{\vartheta})\|_{L^2}
<\infty\quad\hbox{with }\
\langle t\rangle= \sqrt{1+t^2}, \end{eqnarray}
when  the  initial data
 are  perturbed around the equilibrium state $(\bar{\rho},0,\bar{\vartheta})$.
Similar decay estimates  have been established in the half-space or exterior domains (see, for example, \cite{kagei, kobayashi, ma3}).
 Subsequently, the $L^2$ decay estimate \eqref{R-E7} was  improved by  Ponce \cite{ponce} to the $L^p$ one  at the  rate ${\mathcal O}(\langle t\rangle^{-\frac{d}{2}(1-\frac{1}{p})})$ with $2\leq p\leq\infty$ and $\, d=2,3$.
 Inspired by  a series of work   \cite{xujiang2017arma, xujiang2019jmfm, ma2, xujiang2019arxiv, xujiang2019cmp}, the present study  aims at showing  optimal time-decay estimates for \eqref{mm}
within the critical regularity framework of the solutions constructed in \cite{helingbing}.


%
Let  $(\bar{\rho}, \bar u=0, \bar{\theta})$ be an equilibrium state with constants  $\bar{\rho}>0$ and  $\bar{\theta}>0$. We look for the solutions of \eqref{mm} such that
$$(\rho(x, t), u(x,t), \theta(x,t))\to(\bar{\rho}, 0, \bar{\theta}) \quad \mathrm{ for }
\quad|x|\to \infty, \  t>0.$$
For convenience, let  $\bar{\rho}=\bar{\theta}=1$ and define
 $$\mathcal{A}u=\mu\Delta u+(\lambda+\mu)\nabla\div u,\quad  \rho=1+a, \quad \theta=1+\vartheta, \quad I(a)=\frac{a}{1+a}, \quad J(a)=\ln({1+a} ).$$
Then, we can rewrite \eqref{mm} into the following system:
\begin{eqnarray}\label{m}
\left\{\begin{aligned}
&\partial_t a + \div u = - \div(a u)\, ,\\
 &     \partial_t u-\mathcal{A}  u +\nabla (a+\vartheta)= - u\cdot\nabla u-I(a)\mathcal{A}  u+ I(a)\nabla a-\vartheta\nabla J(a),\\
      &\partial_t\vartheta-\kappa\Delta\vartheta+ \div u=- \div(\vartheta u)-\kappa I(a)\Delta \vartheta+\frac{2\mu |D(u)|^2+\lambda(\div u)^2}{1+a},\\
    &(a,u,\vartheta)|_{t=0}=(a_0,u_0,\vartheta_0).
\end{aligned}\right.
\end{eqnarray}

\noindent
For $z \in\mathcal{S}'(\R^3)$,
the low  and high  frequency parts are expressed as
\begin{align}\label{eq:lhf0}
z^\ell\stackrel{\mathrm{def}}{=}\sum_{2^j\leq j_0}\ddj z\quad\hbox{and}\quad
z^h\stackrel{\mathrm{def}}{=}\sum_{2^   j>j_0}\ddj z
\end{align}
for a large  integer $j_0\ge 1.$  
The corresponding  truncated semi-norms are defined  as follows:
\begin{align*}\|z\|^{\ell}_{\dot B^{s}_{p,1}}
\stackrel{\mathrm{def}}{=}  \|z^{\ell}\|_{\dot B^{s}_{p,1}}
\ \hbox{ and }\   \|z\|^{h}_{\dot B^{s}_{p,1}}
\stackrel{\mathrm{def}}{=}  \|z^{h}\|_{\dot B^{s}_{p,1}}.
\end{align*}
The definitions of the operator $\ddj$ and the space  $\dot B^{s}_{p,1}$ are given in the next section.
Without loss of generality, we take $\mu=\lambda=\kappa=1 $ in \eqref{m}.

The present  analysis is closely related to  the following results.
 \begin{theorem}\label{can} ( Danchin and  He \cite{helingbing})
 Let  $2\le p<3$, $(a_0^\ell, u_0^\ell,\ta_0^\ell)\in \dot{B}_{2,1}^{\frac 12}(\R^3)$, $u_0^h\in \dot{B}_{p,1}^{\frac 3p-1}(\R^3)$,  $a_0^h\in \dot{B}_{p,1}^{\frac 3p}(\R^3)$, $\ta_0^h\in \dot{B}_{p,1}^{\frac 3p-2}(\R^3)$. The initial data satisfy the assumption
 \begin{align}\label{new1}
\z_{0} \stackrel{\mathrm{def}}{=}\|(a_0,u_0,\vartheta_0)\|^{\ell}_{\dot B^{\frac 12}_{2,1}}+\|a_0\|^h_{\dot B^{\frac 3p}_{p,1}}+\|u_0\|^h_{\dot B^{\frac 3p-1}_{p,1}}+\|\vartheta_0\|^h_{\dot B^{\frac 3p-2}_{p,1}}\leq \varepsilon
\end{align}
for a  constant $\varepsilon>0$ sufficiently small.
Then
the system \eqref{m} has a unique global solution $(a, u,\ta)$ so that
\begin{align*}
&(a^\ell,u^\ell,\ta^\ell)\in C([0,\infty );{\dot{B}}_{2,1}^{\frac {1}{2}}(\R^3))\cap L^{1}
(0, \infty;{\dot{B}}_{2,1}^{\frac 52}(\R^3)),\\
&u^h\in C([0,\infty );{\dot{B}}_{p,1}^{\frac {3}{p}-1}(\R^3))\cap L^{1}
(0, \infty;{\dot{B}}_{p,1}^{\frac 3p+1}(\R^3)),\\
&a^h\in C([0,\infty );{\dot{B}}_{p,1}^{\frac {3}{p}}(\R^3))\cap L^{1}
(0, \infty;{\dot{B}}_{p,1}^{\frac 3p}(\R^3)),\\
&\ta^h\in C([0,\infty );{\dot{B}}_{p,1}^{\frac {3}{p}-2}(\R^3))\cap L^{1}
(0, \infty;{\dot{B}}_{p,1}^{\frac 3p}(\R^3)).
\end{align*}
Moreover,  there holds
\begin{equation}\label{yaya}
\z(t)+\int^t_0\y(s)\,ds\le C \z_{0},\quad\hbox{for $ t>0$},
\end{equation}
where, $C$ is a generic constant,
\begin{align}
\z(t)\stackrel{\mathrm{def}}{=}&\|(a,u,\ta)\|^\ell_{\dot{B}_{2,1}^\frac12}+\|u\|^h_{\dot{B}_{p,1}^{\frac3p-1}}
+\|a\|^h_{\dot{B}_{p,1}^\frac3p}+\|\ta\|^h_{\dot{B}_{p,1}^{\frac3p-2}},\label{new2}\\
\y(t)\stackrel{\mathrm{def}}{=}&\|(a,u,\vartheta)\|^{\ell}_{\dot{B}^{\frac 52}_{2,1}}
+\|(a,\vartheta)\|^{h}_{\dot{B}^{\frac 3p}_{p,1}}+\|u\|^{h}_{\dot{B}^{\frac 3p+1}_{p,1}}\label{new11}.
\end{align}
\end{theorem}

\begin{theorem}\label{rty} (Danchin and Xu \cite{xujiang2019jmfm})
 In addition to  the conditions of Theorem \ref{can},
 assume that there exists  a small constant $\varepsilon>0$ such that
\begin{equation}\label{tiankgn}
 \|(a_0, u_0, \vartheta_0)\|_{\dot{B}^{\sigma}_{2,1}}^\ell \leq \varepsilon,
\end{equation}
for
\begin{equation}\label{tianmi}
\frac32-\frac6p<\sigma\le -\frac12.
\end{equation}
Then,  for $\ t\geq0$ and $\sigma \le s \le\frac 52$,
\begin{align*}
&\langle t\rangle^{\frac {s-\sigma}2}\|(a,u,\vartheta)(t)\|_{\dot B^s_{2,1}}^\ell
+\langle t\rangle^{2-\sigma}\| a(t)\|_{ \dot B^{\frac 3p}_{p,1}}^h
\nonumber\\
&\quad+\langle t\rangle^{2-\sigma}\| u(t)\|_{ B^{\frac 3p-1}_{p,1}}^h
+\langle t\rangle^{2-\sigma} \|\vartheta(t)\|_{\dot B^{\frac 3p-2}_{p,1}}^h
+t^{2-\sigma}\|(\nabla u,\vartheta)(t)\|_{\dot B^{\frac 3p}_{p,1}}^h\le C\z_{0}.
\end{align*}

\end{theorem}

 The present paper can be regarded as a further study on the previous  theorem without  the additional smallness assumption (\ref{tiankgn}).
  For the compressible Navier-Stokes equations without  (\ref{tiankgn}) (let $\vartheta_0=0$), the optimal time decay has been obtained by Xin and Xu \cite{xujiang2019arxiv}.
Now we  enlarge the range of $\sigma$  in \eqref{tianmi} and remove  \eqref{tiankgn} for the  compressible Navier-Stokes  equations involving  heat-conductive gases. Compared to the compressible Navier-Stokes equations discussed in \cite{xujiang2019arxiv}, the equations \eqref{m} contain the  stronger nonlinear terms $I(a)\Delta \vartheta$ and $I(a)|D(u)|^2$, which lead to analysis difficulties. 

The main result of the present paper reads as follows.
\begin{theorem}\label{th2} In addition to the conditions of Theorem \ref{can},
assume  that   $$ 2 \le p<3,\ \ \frac32-\frac6p< \sigma<\frac 12, \ \ \frac 3p-\frac 32+\sigma\le\alpha\le\frac 3p-1,\ \ \quad\frac 3p-\frac 32+\sigma\le\beta\le\frac 3p-2,$$
the initial data $(a_0^\ell,u_0^\ell,\ta_0^\ell)\in{\dot{B}_{2,1}^{\sigma}}(\R^3)$.
Then the decay rate estimates
\begin{align*}
\|\Lambda^{\alpha} (a,u)\|_{L^p}
\le C(1+t)^{-\frac{3}{2}(\frac 12-\frac 1p)-\frac{\alpha-\sigma}{2}} \ \mbox{ and } \
\|\Lambda^{\beta} \ta\|_{L^p}
\le C(1+t)^{-\frac{3}{2}(\frac 12-\frac 1p)-\frac{\beta-\sigma}{2}}
\end{align*}
hold true for a positive constant $C$ and  the positive differential operator $\Lambda^{s}=\mathcal{F}^{-1}\left|\xi\right|^{s}\mathcal{F}$ with $\mathcal{F}$ the Fourier transformation.
\end{theorem}
\begin{remark}
Compared  with  Theorem \ref{rty} on the decay rates, Theorem \ref{th2}  does not need  the smallness assumption \eqref{tiankgn}. Moreover, the  $\sigma$ range  is larger than that in \eqref{tianmi} and thus  more decay properties can be derived.
\end{remark}

\section{Preliminaries}
For convenience, we use the symbols $\|(a,b)\|_{X}=\|a\|_{X}+\|b\|_{X}$  and $F \lesssim G$, which  represents the inequality $F\leq CG$ for a generic constant $C$.

In this section, we recall some basic facts on Littlewood-Paley theory (see \cite{bcd} for instance).
Let $\chi\ge 0$ and $\varphi\ge 0$ be two smooth radial functions so that
the support of $\chi$ is contained in the ball $\{\xi\in\R^3: |\xi|\leq\frac {4}{3}\}$,   the support of $\varphi$ is contained in  the annulus $\{\xi\in\R^3: \frac {3}{4}\leq|\xi|\leq\frac {8}{3}\}$ and
 \begin{align*}
 \sum_{j \in \mathbb{Z}} \varphi(2^{-j}\xi)=1, \hspace{0.5cm}\forall \xi \neq0.
 \end{align*}
   Let $\mathcal F$ be the Fourier transform.  
The homogeneous dyadic blocks $\dot{\Delta}_{j}$ and the  low-frequency cutoff operators
$\dot{S}_{j}$ are defined for all $j\in\mathbb{Z}$ by
\begin{align*}
\dot{\Delta}_{j}u={\mathcal F}^{-1}(\varphi(2^{-j}\cdot) {\mathcal F}u),\quad\quad
\dot{S}_{j}u={\mathcal F}^{-1}(\chi(2^{-j}\cdot) {\mathcal F}u) .
\end{align*}
Let us remark that, for any homogeneous function $A$ of order 0 smooth outside 0, we have
\begin{equation*}\label{}
\forall p\in[1,\infty],\quad\quad\|\ddj (A(D) f )\|_{L^p}\le C\|\ddj f\|_{L^p}.
\end{equation*}
Denote by $\mathscr{S}_h^{'}(\R^3)$ the space of tempered distributions subject to the condition
$$
\lim_{j\rightarrow -\infty}\dot{S}_j u=0.
$$
Then we have the  decomposition
$$
u = \sum_{j\in \Z} \dot \Delta_j u \ \ \  \forall u\in\mathscr{S}_h^{'}(\mathbb{R}^3).
$$

We recall the definition of homogeneous Besov spaces.
 \begin{definition}
 For  $1 \le p,\ r \le \infty$ and  $-\infty <s <\infty$,
 we set
 \begin{align*}\dot{B}_{p,r}^s(\R^3)=\Big\{u\in \mathscr{S}_h^{'}(\R^3)\ \big| \ \|u\|_{\dot{B}_{p,r}^s}<\infty\Big\}, \end{align*}
where
$$\|u\|_{\dot{B}_{p,r}^s}=\Big(\sum_{j\in \mathbb{Z}} \big(2^{js}\|\dot{\Delta}_ju\|_{L^p}\big)^r\Big)^\frac1r.$$
\end{definition}

Let us now state some Besov space  properties  to be used repeatedly  in this paper.
\begin{lemma}\label{qianrudingli}
(a) For  $s\in \R$ and  $1\le p\le \infty$,  then there holds the embedding
\begin{align*}\|u\|_{\dot{B}_{p,1}^{s}}\lesssim  \|\nabla u\|_{\dot{B}_{p,1}^{s-1}}\lesssim \|u\|_{\dot{B}_{p,1}^{s}} \ \ \mbox{ for } \ \ u\in\dot{B}^{s}_{p,1}(\R^3).
\end{align*}

(b) For   $1\le p\le \infty$ and $s_1,\ s_2\in \R$  with $s_1>s_2$,    then there holds
\begin{align*}
\|u^\ell\|_{\dot{B}^{s_1}_{p,1}}\lesssim\| u^\ell\|_{\dot{B}^{s_2}_{p,1}}, \quad \|u^h\|_{\dot{B}^{s_2}_{p,1}}\lesssim\| u^h\|_{\dot{B}^{s_1}_{p,1}} \  \ \mbox{ for } \ \ u\in\dot{B}^{s_1}_{p,1}(\R^3)\cap\dot{B}^{s_2}_{p,1}(\R^3).
\end{align*}

(c) For  $s\in \R$ and   $1\le p<q\le \infty$, then we have the embedding
\begin{align*}
{\dot{B}_{p,1}^{s}}(\R^3)\hookrightarrow {\dot{B}_{q,1}^{s-\frac{3}{p}-\frac{3}{q}}}(\R^3).
\end{align*}
\end{lemma}

The following Bernstein's lemma will be repeatedly used throughout this paper:
\begin{lemma}\label{bernstein}(see \cite{bcd})
Let $\mathcal{B}$ be a ball and $\mathcal{C}$ an annulus of $\mathbb{R}^3$ centered at the origin. For an integer $0\leq k\leq 2$ and reals
$1\le p \le q\le\infty$, there hold
\begin{align*}
&\sup_{|\alpha|=k}\|\partial^{\alpha}u\|_{L^q}\lesssim \sigma ^{k+3(\frac1p-\frac1q)}\|u\|_{L^p} \ \mbox{ for } \ \ \mathrm{supp} \,{\mathcal F}{u}\subset\sigma  \mathcal{B},\\
&\sigma ^k\|u\|_{L^p}\lesssim \sup_{|\alpha|=k}\|\partial^{\alpha}u\|_{L^p}
\lesssim \sigma ^{k}\|u\|_{L^p} \ \mbox{ for } \ \ \mathrm{supp} \,{\mathcal F}{u}\subset\sigma  \mathcal{C}
 \end{align*} with respect to scaling parameter $\sigma>0$.
\end{lemma}

\begin{lemma}\label{neicha} (\cite{bcd})
For $0<s_1<s_2,$ $0\leq \theta\leq 1$ and  $1\leq p\leq\infty$, there hold the following interpolation inequality
\begin{align*}
\|u\|_{\dot{B}_{p,1}^s}\lesssim\|u\|^\theta_{\dot{B}_{p,1}^{s_1}}
\|u\|^{1-\theta}_{\dot{B}_{p,1}^{s_2}}\quad \mbox{ with } \,\, s=\theta s_1+(1-\theta)s_2.
\end{align*}

\end{lemma}
The Bony decomposition  is
very effective in the estimate of nonlinear terms in fluid motion equations.
Here, we recall the decomposition in the homogeneous context:
\begin{align*}
uv=\dot{T}_uv+\dot{T}_vu+\dot{R}(u,v),
\end{align*}
where
$$\dot{T}_uv\stackrel{\mathrm{def}}{=}\sum_{j\in \Z}\dot{S}_{j-1}u\dot{\Delta}_jv \ \ \mbox{ and }\ \ \hspace{0.5cm}\dot{R}(u,v)\stackrel{\mathrm{def}}{=}\sum_{j\in \Z}
\dot{\Delta}_ju\widetilde{\dot{\Delta}}_jv \  \ \mbox{ with } \ \ \widetilde{\dot{\Delta}}_jv\stackrel{\mathrm{def}}{=}\sum_{|j-j'|\le1}\dot{\Delta}_{j'}v.$$

The paraproduct $\dot{T}$ and the remainder $\dot{R}$ operators satisfy the following
continuous properties.

\begin{lemma}(see \cite{bcd})\label{fangji}
For  $s, t \in\mathbb{R}$,  $1\leq p_1, p_2\leq\infty$ and $\frac{1}{p}=\frac{1}{p_1}+\frac{1}{p_2}$,  there hold the inequalities
\begin{align*}
&\|\dot T_uv\|_{\dot{B}_{p,1}^{s+t}}\lesssim \|u\|_{\dot{B}_{p_1,1}^{s}}\|v\|_{ \dot{B}_{p_2,1}^t},\,\,\, s\le 0,
\\
&\|\dot R(u,v)\|_{\dot{B}_{p,1}^{s+t}}\lesssim \|u\|_{\dot{B}_{p_1,1}^{s}}\|v\|_{ \dot{B}_{p_2,1}^t},\,\,\, s+t>0,\\
&\|\dot R(u,v)\|_{\dot{B}_{p,\infty}^{0}}\lesssim \|u\|_{\dot{B}_{p_1,1}^{s}}\|v\|_{ \dot{B}_{p_2,\infty}^{-s}}.
 \end{align*}

\end{lemma}
The following product law  plays central roles when we estimate the couple terms appeared in the equations:
\begin{lemma}\label{daishu}( \cite[Proposition A.1]{xujiang2017arma})
 Assume that $1\leq p,\, q\leq \infty$,\,
 $$s_1\leq \frac {3}{q}, \ \ s_2\leq 3\min\left\{\frac 1p\, ,\,\frac 1q\right\} \mbox{ and } s_1+s_2>3\max\left\{0,\frac 1p +\frac 1q -1\right\}.$$ Then we have,  for $ (f,g)\in\dot{B}_{q,1}^{s_1}({\mathbb R} ^3)\times\dot{B}_{p,1}^{s_2}({\mathbb R} ^3)$,
\begin{align*}
\|fg\|_{\dot{B}_{p,1}^{s_1+s_2 -\frac {3}{q}}}\lesssim \|f\|_{\dot{B}_{q,1}^{s_1}}\|g\|_{\dot{B}_{p,1}^{s_2}}.
\end{align*}
\end{lemma}
In order to deal with composition functions in the Besov spaces, we also need the following proposition:
\begin{proposition}\label{pro}
(see \cite{bcd}) For $ \alpha<0<\beta$, let  $G$ be  a smooth function defined on the  open interval $(\alpha,\beta)$ so that   $G(0)=0.$
Assume that  $f\!: \R\mapsto {\mathcal I} \subset (\alpha,\beta)$ for an interval ${\mathcal I}$. Then we have    the estimate
\begin{align}\label{new44}
\|G(f)\|_{\dot B^{s}_{p,1}}\lesssim\|f\|_{\dot B^s_{p,1}}\,\,\,\mbox{ for }\,\,\, 1\leq p\leq\infty,\,\,\,\, s>0.
\end{align}
\end{proposition}

\section{The proof of Theorem \ref{th2}}

The proof is to show the decay estimate of the solution $(a,u,\vartheta)$ given in Theorem \ref{can}. This decay estimate is simply derived  from  a Lyapunov-type differential inequality for the solution. This inequality lies on the viability of the uniform bound
\begin{equation}\label{new3}\|(a,u,\vartheta)\|^{\ell}_{\dot{B}^{\sigma}_{2,1}}\le C \,\,\mbox{ for }\,\, \sigma<\frac12
\end{equation}
This bound can be derived from an energy estimate of (\ref{m}). To show the energy estimate, we starts with the following energy estimate the linearized equation system of (\ref{m}).

\subsection{Low frequency energy estimate of the linearized system of \eqref{m}}
For  $(a,u,\vartheta)$  the smooth solution to the following linearized  system of (\ref{m}):
\begin{eqnarray}\label{lm}
\left\{\begin{aligned}
&\partial_t a + \div u = f_1\, ,\\
 &     \partial_t u-\Delta  u +\nabla (a+\vartheta)= f_2,\\
      &\partial_t\vartheta-\Delta\vartheta+ \div u=f_3,\\
    &(a,u,\vartheta)|_{t=0}=(a_0,u_0,\vartheta_0),
\end{aligned}\right.
\end{eqnarray}
we show the energy estimate,  for  $\gamma\in\R$,
\begin{align}\label{qi1}
\|(a,u,\vartheta)\|^{\ell}_{\dot{B}^{\gamma}_{2,1}}+\int_0^t\|(a,u,\vartheta)\|^{\ell}_{\dot{B}^{\gamma+2}_{2,1}}\,ds
\lesssim\|(a_0,u_0,\vartheta_0)\|^{\ell}_{\dot{B}^{\gamma}_{2,1}}+\int_0^t\|(f_1,f_2,f_3)\|^{\ell}_{\dot{B}^{\gamma}_{2,1}}\,ds.
\end{align}

Indeed,  let $\omega=\Lambda^{-1}\mathrm{curl} \,u$ be the incompressible part of $u$ and $v=\Lambda^{-1}\div u$  be the compressible part of $u$. We see that $\omega$ satisfies the heat equation:
\begin{equation}\label{qi2}
\partial_{t}\omega-\Delta\omega=\Lambda^{-1}\mathrm{curl}f_2,\quad\omega(0)=\omega_0.
\end{equation}
A standard energy argument applied to \eqref{qi2} implies
\begin{align}\label{qi3}
\|\omega\|^{\ell}_{\dot{B}^{\gamma}_{2,1}}+\int_0^t\|\omega\|^{\ell}_{\dot{B}^{\gamma+2}_{2,1}}\,ds
\lesssim\|\omega_0\|^{\ell}_{\dot{B}^{\gamma}_{2,1}}+\int_0^t\|f_2\|^{\ell}_{\dot{B}^{\gamma}_{2,1}}\,ds.
\end{align}
On the other hand, it is easy to check that $(\ddj a,\ddj v,\ddj \ta)$ satisfies the equations
\begin{eqnarray}\label{qiqi}
\left\{\begin{aligned}
&\partial_t \ddj a + \La \ddj v = \ddj f_1\, ,\\
 &     \partial_t \ddj v-\Delta \ddj  v -\ddj \La(a+\vartheta)=\La^{-1}\ddj \div f_2,\\
      &\partial_t\ddj \vartheta-\Delta\ddj \vartheta+ \La\ddj v=\ddj f_3.
\end{aligned}\right.
\end{eqnarray}

Taking the $L^2$ inner product of \eqref{qiqi} with $(\ddj a, \ddj v, \ddj\ta)$,  and using the following cancellation
$$\int_{\R^3}\ddj \La v\cdot \ddj a\,dx+\int_{\R^3}\ddj \La v\cdot \ddj \ta\,dx-\int_{\R^3}\ddj\La(a+\vartheta)\cdot \ddj v\,dx=0,$$
we  get
\begin{align}\label{qi4}
&\frac 12 \frac{d}{dt}(\|\ddj a\|^2_{L^2}+\|\ddj v\|^2_{L^2}+\|\ddj \ta\|^2_{L^2} )+\|\ddj \Lambda v\|^2_{L^2}+\|\ddj \Lambda \ta\|^2_{L^2}\nonumber\\
&\quad=\int_{\R^3}\ddj f_1\cdot \ddj a\,dx+\int_{\R^3}\ddj\La^{-1} \div f_2\cdot \ddj v\,dx+\int_{\R^3}\ddj f_3\cdot \ddj \ta\,dx.
\end{align}
Applying   $\La$ to the first equation of \eqref{qiqi} gives
\begin{align}\label{qi5}
\partial_t \ddj\La  a -\ddj\Delta v = \ddj\La f_1.
\end{align}
Taking the $L^2$ inner product of \eqref{qi5} with $\ddj \La a$  implies
\begin{align}\label{qi6}
&\frac 12 \frac{d}{dt}\|\ddj \La a\|^2_{L^2}-\int_{\R^3}\ddj \Delta v\cdot \ddj \La a\,dx=\int_{\R^3}\ddj \La f_1\cdot \ddj \La a\,dx.
\end{align}

To find the hidden dissipation of $a$ in the low frequency, we get  by
testing  the second equation of  \eqref{qiqi}  by $\ddj\La a $ and \eqref{qi5} by $\ddj v$ respectively that
\begin{align}\label{qi7}
& \frac{d}{dt}\left(\int_{\R^3}\ddj  v\cdot \ddj \La a\,dx\right)+\|\ddj \Lambda v\|^2_{L^2}-\|\ddj \Lambda a\|^2_{L^2}-\int_{\R^3}\ddj \La\ta \cdot \ddj \La a\,dx-\int_{\R^3}\ddj \Delta v\cdot \ddj \La a\,dx
\nonumber\\
&\quad=\int_{\R^3}\ddj \La f_1\cdot \ddj v\,dx+\int_{\R^3}\ddj\La^{-1} \div f_2\cdot \ddj \La a\,dx.
\end{align}
%
Multiplying (\ref{qi4}) by $2$ and \eqref{qi7} by $-1$ respectively, and summing up the resultant equations with  \eqref{qi6},  we have the energy equality involving the dissipation for $a$:
\begin{align}\label{qi10}
&\frac 12 \frac{d}{dt}\mathcal{L}_{j}^2+\|\ddj \Lambda v\|^2_{L^2}+2\|\ddj \Lambda \ta\|^2_{L^2}+\|\ddj \Lambda a\|^2_{L^2}+\int_{\R^3}\ddj \La\ta \cdot \ddj \La a\,dx\nonumber\\
&\quad=2\int_{\R^3}\ddj f_1\cdot \ddj a\,dx
+2\int_{\R^3}\ddj\La^{-1} \div f_2\cdot \ddj v\,dx+2\int_{\R^3}\ddj f_3\cdot \ddj \ta\,dx\nonumber
\\ &-\int_{\R^3}\ddj \La f_1\cdot \ddj v\,dx-\int_{\R^3}\ddj\La^{-1} \div f_2\cdot \ddj \La a\,dx
+\int_{\R^3}\ddj \La f_1\cdot \ddj \La a\,dx,
\end{align}
where
$$\mathcal{L}_{j}\stackrel{\mathrm{def}}{=} \left(2\|\ddj a\|^2_{L^2}+\|\ddj\La a\|^2_{L^2}+2\|\ddj v\|^2_{L^2}+2\|\ddj\ta\|^2_{L^2}-2\int_{\R^3}\ddj  v\cdot \ddj \La a\,dx\right)^\frac12.$$
It readily seen that the equivalence
\begin{align}
\mathcal{L}^2_{j}\lesssim  \|\ddj a\|^2_{L^2}+\|\ddj v\|^2_{L^2}+\|\ddj\ta\|^2_{L^2}  \lesssim \mathcal{L}^2_{j}
\label{qi11}
\end{align}
holds true for low frequency   components ($j\leq j_0$).
Hence, by Young inequality,  H\"older inequality and Bernstein inequalities, we have
\begin{equation}\label{qi12}
\frac 12 \frac{d}{dt}\mathcal{L}^2_{j}+2^{2j}\mathcal{L}^2_{j}\lesssim \|(\ddj f_1,\ddj f_2,\ddj f_3)\|_{L^2}\mathcal{L}_{j}, \,\, j<j_0,
\end{equation} which implies that
\begin{align} \label{qi13}
&\|(\ddj a,\ddj v,\ddj\ta)\|_{L^2}+2^{2j}\int_0^t\|(\ddj a,\ddj v,\ddj\ta)\|_{L^2}\,ds\nonumber\\
&\quad\lesssim\|(\ddj a_0,\ddj v_0,\ddj\ta_0)\|_{L^2}+\int_0^t\|(\ddj f_1,\ddj f_2,\ddj f_3)\|_{L^2}\,ds, \quad \hbox{ $j\leq j_0$.}
\end{align}
 Multiplying  \eqref{qi13} by $2^{\gamma j}$ and summing up the resultant inequalities with respect to  $j\le j_0$, we have
\begin{align}\label{qi14}
\|(a,v,\vartheta)\|^{\ell}_{\dot{B}^{\gamma}_{2,1}}+\int_0^t\|(a,v,\vartheta)\|^{\ell}_{\dot{B}^{\gamma+2}_{2,1}}\,ds
\lesssim\|(a_0,v_0,\vartheta_0)\|^{\ell}_{\dot{B}^{\gamma}_{2,1}}+\int_0^t\|(f_1,f_2,f_3)\|^{\ell}_{\dot{B}^{\gamma}_{2,1}}\,ds
\end{align}
which combined with  \eqref{qi3} and \eqref{qi14} gives the desired  the low frequencies estimate \eqref{qi1}.

\subsection{Nonlinear estimates for showing the uniform boundedness \eqref{new3}}

As
$(a,u,\vartheta)$  is  the global small solution of \eqref{m} given  in Theorem \ref{can} and \eqref{lm} is the linearized equation system of \eqref{m}, the application of (\ref{qi1}) to \eqref{m} with $\gamma=\sigma$ for   $\frac32-\frac6p< \sigma<\frac 12$ gives
\begin{align}
&\|(a,u,\vartheta)\|^{\ell}_{\dot{B}^{\sigma}_{2,1}}+\int_0^t\|(a,u,\vartheta)\|^{\ell}_{\dot{B}^{\sigma+2}_{2,1}}\,ds\nonumber\\
&\quad\lesssim\|(a_0,u_0,\vartheta_0)\|^{\ell}_{\dot{B}^{\sigma}_{2,1}}+\int_0^t\|u\cdot \nabla u\|^{\ell}_{\dot{B}^{\sigma}_{2,1}}\,ds+\int_0^t\|\div(\vartheta u)\|^{\ell}_{\dot{B}^{\sigma}_{2,1}}\,ds+\int_0^t\|\div (au)\|^{\ell}_{\dot{B}^{\sigma}_{2,1}}\,ds
\nonumber\\
&\quad\quad+\int_0^t\|I(a)\mathcal{A} u\|^{\ell}_{\dot{B}^{\sigma}_{2,1}}\,ds+\int_0^t\|I(a)\nabla a\|^{\ell}_{\dot{B}^{\sigma}_{2,1}}\,ds+\int_0^t\|\vartheta\nabla J(a)\|^{\ell}_{\dot{B}^{\sigma}_{2,1}}\,ds
\nonumber\\
&\quad\quad+\int_0^t\| I(a)\Delta\vartheta\|^{\ell}_{\dot{B}^{\sigma}_{2,1}}\,ds+\int_0^t\big\|(1+I(a))\big(2|Du|^2+(\div u)^2\big)\big\|^{\ell}_{\dot{B}^{\sigma}_{2,1}}\,ds.\label{jiang2}
\end{align}

To deal with the nonlinear terms on the right-hand side of the previous equation, we need the following  product laws:
\begin{align}
\|fg\|_{\dot{B}_{2,1}^{\sigma}}\lesssim &\|f\|_{\dot{B}_{2,1}^{\sigma}}\|g\|_{\dot{B}_{p,1}^{\frac 3p}},\quad \hbox{$ -\frac 3p<\sigma\le \frac 3p,\quad 2\le p<3$},\label{ci1}\\
\|f g^h\|_{\dot B^{\sigma}_{2,1}}^\ell \lesssim&\|f\|_{\dot B^{\frac 3p-1}_{p,1}}\|g^h\|_{\dot B^{\frac 3p-1}_{p,1}},\quad \hbox{$ \frac32-\frac6p< \sigma\le \frac 3p,\quad 2\le p<3$}.\label{ci2}
\end{align}
The first one is given by  Lemma \ref{daishu}. To prove the second one, we use   Bony's decomposition:
$$
fg^h= \dot{T}_{g^h}f+\dot{R}(g^h,f)+\dot{T}_{f}g^h.
$$
 Applying Lemma \ref{fangji} and the condition $1-\frac 3p<0$ implies that
\begin{align}\label{jiang998}
\|\dot{T}_{g^h}f\|_{\dot B^{\sigma}_{2,1}}^\ell \lesssim&\|\dot{T}_{g^h}f\|_{\dot B^{\frac32-\frac6p}_{2,1}}^\ell \lesssim \|g^h\|_{\dot B^{1-\frac 3p}_{p,1}}\|f\|_{\dot B^{\frac 3p-1}_{p,1}}
\lesssim \|g^h\|_{\dot B^{\frac 3p-1}_{p,1}}\|f\|_{\dot B^{\frac 3p-1}_{p,1}}
\end{align}
where we have used the high  frequency property of $g^h$ and the  fact $1-\frac 3p<\frac 3p-1$ in the last inequality.
 Similarly,
by using the low  frequency property and the condition $\frac32-\frac6p< \sigma<\frac 12$, the rest term can be estimated from Lemma \ref{fangji} that
\begin{align}\label{jiang99899}
\|\dot{R}(g^h,f)\|^\ell_{\dot B^{\sigma}_{2,1}}
\lesssim&\|\dot{R}(g^h,f)\|_{\dot B^{\frac32-\frac6p}_{2,\infty}}\lesssim\|\dot{R}(g^h,f)\|_{\dot B^{0}_{p/2,\infty}}\nonumber\\
\lesssim& \|g^h\|_{\dot B^{1-\frac 3p}_{p,1}}\|f\|_{\dot B^{\frac 3p-1}_{p,1}}
\lesssim \|g^h\|_{\dot B^{\frac 3p-1}_{p,1}}\|f\|_{\dot B^{\frac 3p-1}_{p,1}}.
\end{align}
Moreover, with the aid of the  H\"{o}lder inequality and Bernstein's inequality, we have
\begin{align}
\|\dot{T}_fg^h\|^\ell_{\dot{B}^{\sigma}_{2,1}}&\lesssim \sum_{j \le j_0}\sum_{|j-k|\leq 4}2^{j(\frac32-\frac6p)}\|\dot{\Delta}_j\big(\dot{S}_{k-1}f\dot{\Delta}_{k}g^h\big)\|_{L^2}
\nonumber\\&\lesssim \sum_{j \le j_0}\sum_{|j-k|\leq 4}2^{j(\frac32-\frac6p)}\|\dot{S}_{k-1}f\|_{L^{\frac{2p}{p-2}}}\|\dot{\Delta}_{k}g^h\|_{L^p}
\nonumber\\&\lesssim\sum_{j \le j_0}\sum_{|j-k|\leq 4}2^{j(\frac32-\frac6p)}(\sum_{k'\leq k-2}\|\dot{\Delta}_{k'} f\|_{L^{{\frac{2p}{p-2}}}})\|\dot{\Delta}_kg^h\|_{L^p}\nonumber
\end{align}
which together with
a similar derivation of \eqref{jiang998}
yields
\begin{align}\label{jiang999}
\|\dot{T}_fg^h\|^\ell_{\dot{B}^{\sigma}_{2,1}}&\lesssim\sum_{j \le j_0}\sum_{|j-k|\leq 4}2^{j(\frac32-\frac6p)}(\sum_{k'\leq k-2}2^{(\frac{6}{p}-\frac{3}{2})k'}\|\dot{\Delta}_{k'} f\|_{L^{p}})\|\dot{\Delta}_kg^h\|_{L^p}
\nonumber\\&\lesssim \|f\|_{\dot{B}^{\frac 3p-1}_{p,1}}\sum_{j \le j_0}\sum_{|j-k|\leq 4}2^{(j-k)(\frac 32-\frac 6p)}2^{k{(1-\frac 3p)}}\|\dot{\Delta}_{k}g^h\|_{L^p}
\nonumber\\&\lesssim \|f\|_{\dot{B}^{\frac 3p-1}_{p,1}}\|g^h\|_{\dot{B}^{1-\frac 3p}_{p,1}}\lesssim \|f\|_{\dot{B}^{\frac 3p-1}_{p,1}}\|g^h\|_{\dot{B}^{\frac 3p-1}_{p,1}}.
\end{align}
The combination of \eqref{jiang998}, \eqref{jiang99899} and \eqref{jiang999}  yields \eqref{ci2}.

 In  estimating the nonlinear items of \eqref{jiang2} by using \eqref{ci1}, \eqref{ci2} and the decomposition $u=u^\ell+u^h$, we find that
\begin{align} 
\|u\cdot \nabla u\|^{\ell}_{\dot{B}^{\sigma}_{2,1}}
\lesssim&\|u^\ell\cdot \nabla u^\ell\|^{\ell}_{\dot{B}^{\sigma}_{2,1}}+\|u^h\cdot \nabla u^\ell\|^{\ell}_{\dot{B}^{\sigma}_{2,1}} +\|u\cdot \nabla u^h\|^{\ell}_{\dot{B}^{\sigma}_{2,1}}
\nonumber\\
\lesssim&(\|\nabla u^{\ell}\|_{\dot{B}^{\frac3p}_{p,1}}+\| u^h\|_{\dot{B}^{\frac3p}_{p,1}})\| \nabla u^\ell\|_{\dot{B}^{\sigma}_{2,1}} +(\|u^\ell\|_{\dot B^{\frac 3p-1}_{p,1}}+\|u^h\|_{\dot B^{\frac 3p-1}_{p,1}})\|\nabla u^h\|_{\dot B^{\frac 3p-1}_{p,1}}.
\nonumber
\end{align}
By using Bernstein's estimate and the properties of low and high frequencies, the previous estimate becomes
\begin{align}\label{jiang4}
\|u\cdot \nabla u\|^{\ell}_{\dot{B}^{\sigma}_{2,1}}
\lesssim& (\|u\|^{\ell}_{\dot{B}^{\frac 52}_{2,1}}+\|u\|^{h}_{\dot{B}^{\frac 3p+1}_{p,1}})\|u\|^{\ell}_{\dot{B}^{\sigma}_{2,1}}+(\|u\|^{\ell}_{\dot B^{\frac 12}_{2,1}}+\|u\|^h_{\dot B^{\frac 3p-1}_{p,1}})\|u\|^h_{\dot B^{\frac 3p+1}_{p,1}}\nonumber\\
\lesssim
&\y(t)\|u\|^{\ell}_{\dot{B}^{\sigma}_{2,1}}+\z(t)\y(t),
\end{align}
for $\y(t)$ and $\z(t)$ defined in \eqref{new11} and \eqref{new1}.

By  \eqref{ci1}, Lemma \ref{bernstein}  and the properties of low and high frequencies, we estimate the following nonlinear terms
\begin{align}\label{jiang6}
&\|u\!\cdot\!\nabla \vartheta^{\ell}\|^{\ell}_{\dot{B}^{\sigma}_{2,1}}\!+\!\|\vartheta\,\mathrm{div}\,u^{\ell}\|^{\ell}_{\dot{B}^{\sigma}_{2,1}}
\nonumber\\&\quad
\lesssim\|u^\ell\|_{\dot{B}^\sigma_{2,1}} \|\nabla \vartheta^\ell\|_{\dot{B}^\frac3p_{p,1}}
\!+\!\|u\|^h_{\dot B^\frac3p_{p,1}}\|\nabla\vartheta^{\ell}\|_{\dot B^\sigma_{2,1}}
\!+\!\|\vartheta^\ell\|_{\dot B^\sigma_{2,1}}\|\mathrm{div}\,u^\ell\|_{\dot B^\frac3p_{p,1}}
\!+\!\|\vartheta^h\|_{\dot B^\frac3p_{p,1}}\|\mathrm{div}\,u^{\ell}\|_{\dot B^\sigma_{2,1}}\nonumber\\
&\quad\lesssim\|\vartheta^{\ell}\|_{\dot{B}^{\frac 52}_{2,1}}\|u^{\ell}\|_{\dot{B}^{\sigma}_{2,1}}\!+\!\|u\|^{h}_{\dot{B}^{\frac 3p\!+\!1}_{p,1}}\|\vartheta^{\ell}\|_{\dot{B}^{\sigma}_{2,1}}\!+\!\|u\|^{\ell}_{\dot{B}^{\frac 52}_{2,1}}\|\vartheta\|^{\ell}_{\dot{B}^{\sigma}_{2,1}}\!+\!\|\vartheta\|^{h}_{\dot{B}^{\frac 3p}_{p,1}}\|u\|^{\ell}_{\dot{B}^{\sigma}_{2,1}}.
\end{align}
By  \eqref{ci2} and Lemma \ref{neicha},
\begin{align}\label{jiang7}
\|u\cdot&\nabla \vartheta^h\|_{\dot B^{\sigma}_{2,1}}^\ell\!+\!\|\vartheta\mathrm{div}u^{h}\|_{\dot B^{\sigma}_{2,1}}^\ell
\nonumber\\
\lesssim& (\|u\|^{\ell}_{\dot B^{\frac3p-1}_{p,1}}\!+\!\|u\|^h_{\dot B^{\frac 3p-1}_{p,1}})\|\nabla \vartheta\|^h_{\dot B^{\frac 3p-1}_{p,1}}
 \!+\!(\|\vartheta\|^\ell_{\dot B^{\frac 3p-1}_{p,1}}\!+\!\|\vartheta\|^h_{\dot B^{\frac 3p-1}_{p,1}})\|\mathrm{div}u\|^h_{\dot B^{\frac 3p-1}_{p,1}}\nonumber
\\
\lesssim& (\|u\|^{\ell}_{\dot B^{\frac 12}_{2,1}}\!+\!\|u\|^h_{\dot B^{\frac 3p-1}_{p,1}})\|\vartheta\|^h_{\dot B^{\frac 3p}_{p,1}}\!+\!\|\vartheta\|^{\ell}_{\dot B^{\frac 12}_{2,1}}\|u\|^h_{\dot B^{\frac 3p+1}_{p,1}}\!+\!\big(\|\vartheta\|^h_{\dot B^{\frac 3p-2}_{p,1}}\|\vartheta\|^h_{\dot B^{\frac 3p}_{p,1}}\|u\|^h_{\dot B^{\frac 3p-1}_{p,1}}\|u\|^h_{\dot B^{\frac 3p+1}_{p,1}}\big)^{\frac12}\nonumber\\
\lesssim& (\|u\|^{\ell}_{\dot B^{\frac 12}_{2,1}}\!+\!\|u\|^h_{\dot B^{\frac 3p-1}_{p,1}})\|\vartheta\|^h_{\dot B^{\frac 3p}_{p,1}}\!+\!\|\vartheta\|^{\ell}_{\dot B^{\frac 12}_{2,1}}\|u\|^h_{\dot B^{\frac 3p+1}_{p,1}}\!+\!\|\vartheta\|^h_{\dot B^{\frac 3p-2}_{p,1}}
\|u\|^h_{\dot B^{\frac 3p+1}_{p,1}}\!+\!
\|\vartheta\|^h_{\dot B^{\frac 3p}_{p,1}}\|u\|^h_{\dot B^{\frac 3p-1}_{p,1}}.
\end{align}
Combining  \eqref{jiang6} and \eqref{jiang7} with the identity
\begin{align}\label{jiang5}
\div (\vartheta u)=u\cdot\nabla \vartheta^\ell+\vartheta\mathrm{div}u^{\ell}+u\cdot\nabla \vartheta^h+\vartheta\mathrm{div}u^{h},
\end{align}
gives
\begin{align}\label{jiang8}
\|\div (\vartheta u)\|^{\ell}_{\dot{B}^{\sigma}_{2,1}}
\lesssim&(\|\vartheta\|^{\ell}_{\dot{B}^{\frac 52}_{2,1}}\!+\!\|\vartheta\|^{h}_{\dot{B}^{\frac 3p}_{p,1}}\!+\!\|u\|^{\ell}_{\dot{B}^{\frac 52}_{2,1}}\!+\!\|u\|^{h}_{\dot{B}^{\frac 3p+1}_{p,1}})\|(\vartheta,u)\|^{\ell}_{\dot{B}^{\sigma}_{2,1}}\nonumber\\
+&(\|\vartheta\|^{\ell}_{\dot{B}^{\frac 12}_{2,1}}\!+\!\|\vartheta\|^{h}_{\dot{B}^{\frac 3p-2}_{p,1}})\|u\|^h_{\dot B^{\frac 3p+1}_{p,1}}
\!+\!(\|u\|^{\ell}_{\dot B^{\frac 12}_{2,1}}\!+\!\|u\|^h_{\dot B^{\frac 3p-1}_{p,1}})\|\vartheta\|^h_{\dot B^{\frac 3p}_{p,1}}\nonumber\\
\lesssim
&\y(t)\|(\vartheta,u)\|^{\ell}_{\dot{B}^{\sigma}_{2,1}}+\y(t)\z(t).
\end{align}

Similarly, we have
\begin{align}\label{jiang9}
\|\div (au)\|^{\ell}_{\dot{B}^{\sigma}_{2,1}}\lesssim
&\y(t)\|(a,u)\|^{\ell}_{\dot{B}^{\sigma}_{2,1}}+\y(t)\z(t).
\end{align}

Moreover, to deal with the nonlinear term $I(a)\nabla a$, we have to use the estimate (\ref{new44}) in a Besov space $\dot B^s_{p,1}$ for $s>0$. However, for the Besov space $\dot B^\sigma_{2,1}$ in the present estimation, the condition  $\sigma>0$ cannot be guaranteed. Thus in order to make use of Proposition \ref{pro} or (\ref{new44}), we  employ  the formulation $I(a)=a-aI(a)$ and  argue    in the same way as  the derivation of \eqref{jiang4}  to obtain that
\begin{align}\label{jiang11}
\|I(a)\nabla a\|^{\ell}_{\dot{B}^{\sigma}_{2,1}}
\lesssim&\| a^{\ell}\|_{\dot{B}^{\sigma}_{2,1}} \|\nabla a^\ell\|_{\dot{B}^{\frac 3p}_{p,1}}
  \!+\!\|a^h\|_{\dot{B}^{\frac 3p}_{p,1}}\|\nabla a^{\ell}\|^{\ell}_{\dot{B}^{\sigma}_{2,1}}
  \nonumber\\
  &\quad\quad+\|aI(a) \|_{\dot{B}^{\frac 3p}_{p,1}}\|a\|^{\ell}_{\dot{B}^{\sigma}_{2,1}} \!+\!\|I(a)\|_{\dot B^{\frac 3p-1}_{p,1}}\|\nabla a\|^h_{\dot B^{\frac 3p-1}_{p,1}}
\nonumber\\ \lesssim&
(\|a\|^{\ell}_{\dot{B}^{\frac 52}_{2,1}}\!+\!\|a\|^{h}_{\dot{B}^{\frac 3p}_{p,1}})\|a\|^{\ell}_{\dot{B}^{\sigma}_{2,1}}\!+\!\|a\|^2_{\dot{B}^{\frac 3p}_{p,1}}\|a\|^{\ell}_{\dot{B}^{\sigma}_{2,1}}
\!+\!\|a\|_{\dot B^{\frac 3p-1}_{p,1}}\|\nabla a\|^h_{\dot B^{\frac 3p-1}_{p,1}}
\nonumber\\ \lesssim&
(\|a\|^{\ell}_{\dot{B}^{\frac 52}_{2,1}}+\|a\|^{h}_{\dot{B}^{\frac 3p}_{p,1}}
+(\|a\|^h_{\dot{B}^{\frac 3p}_{p,1}})^2)\|a\|^{\ell}_{\dot{B}^{\sigma}_{2,1}}+\Big(\|a\|^{\ell}_{\dot B^{\frac 12}_{2,1}}+\|a\|^h_{\dot B^{\frac 3p}_{p,1}}\Big)\|a\|^h_{\dot B^{\frac 3p}_{p,1}}\nonumber
\\ \lesssim &\y(t)(1+\z(t))\|a\|^{\ell}_{\dot{B}^{\sigma}_{2,1}}+\z(t)\y(t).
\end{align}
%
Similarly, we have
\begin{align}\label{jiang14}
 \|I(a)\mathcal{A} u\|_{\dot B^{\sigma}_{2,1}}^\ell
\lesssim&
\big(\|u\|^{\ell}_{\dot{B}^{\frac 52}_{2,1}}+\|a\|^{h}_{\dot{B}^{\frac 3p}_{p,1}}+\|a\|^{\ell}_{\dot{B}^{\frac 12}_{2,1}}\|a\|^{\ell}_{\dot{B}^{\frac 52}_{2,1}}+(\|a\|^h_{\dot{B}^{\frac 3p}_{p,1}})^2\big)\|(a,u)\|^{\ell}_{\dot{B}^{\sigma}_{2,1}}\nonumber\\
\lesssim
&\y(t)\z(t)+\y(t)(1+\z(t))\|(a,u)\|^{\ell}_{\dot{B}^{\sigma}_{2,1}}.
\end{align}

Furthermore, to estimate the nonlinear  term $\vartheta\nabla J(a)$, we continue to use the argument for  the derivation of \eqref{jiang4} and Lemma \ref{neicha} to obtain
\begin{align}
\| \vartheta\nabla J(a)\|^\ell_{\dot B^{\sigma}_{2,1}}
\lesssim&\| \vartheta^\ell\nabla( J(a))^\ell\|_{\dot B^{\sigma}_{2,1}}^\ell+\| \vartheta^\ell\nabla( J(a))^h\|_{\dot B^{\sigma}_{2,1}}^\ell+\|\vartheta^h\nabla J(a)\|_{\dot B^{\sigma}_{2,1}}^\ell\nonumber\\
 \lesssim&
 \| \vartheta^\ell\|_{\dot B^{\sigma}_{2,1}}\|\nabla(J(a))^\ell\|_{\dot B^{\frac 3p}_{p,1}}+\|\vartheta^\ell\|_{\dot B^{\frac 3p-1}_{p,1}}\|\nabla( J(a))^h\|_{\dot B^{\frac 3p-1}_{p,1}}+\|\nabla J(a)\|_{\dot B^{\frac 3p-1}_{p,1}}\|\vartheta^h\|_{\dot B^{\frac 3p-1}_{p,1}}\nonumber\\
  \lesssim&
 \| \vartheta^\ell\|_{\dot B^{\sigma}_{2,1}}\|(J(a))^\ell\|_{\dot B^{\frac 3p+1}_{p,1}}+\|\vartheta^\ell\|_{\dot B^{\frac 3p-1}_{p,1}}\|( J(a))^h\|_{\dot B^{\frac 3p}_{p,1}}+\|a\|_{\dot B^{\frac 3p}_{p,1}}\|\vartheta^h\|_{\dot B^{\frac 3p-1}_{p,1}}\nonumber\\
 \lesssim&
 \| \vartheta^\ell\|_{\dot B^{\sigma}_{2,1}}\|a\|^\ell_{\dot B^{\frac 52}_{2,1}}+\|\vartheta^\ell\|_{\dot B^{\frac 12}_{2,1}}\|a\|^h_{\dot B^{\frac 3p}_{p,1}}
 +(\|a^\ell\|_{\dot B^{\frac 12}_{2,1}}+\|a^h\|_{\dot B^{\frac 3p}_{p,1}})\|\vartheta^h\|_{\dot B^{\frac 3p}_{p,1}}\nonumber\\
  \lesssim& \y(t)\z(t)+\y(t)\|\vartheta\|^{\ell}_{\dot{B}^{\sigma}_{2,1}}.
\label{jiang15}
\end{align}

Next attempt is to consider  the nonlinear  term  $I(a)\Delta\vartheta $ involving composition functions and  is more elaborate.
Its component involving the low frequency $\vartheta^\ell$ is estimated as 
\begin{align}
\| I(a)\Delta\vartheta^\ell\|_{\dot B^{\sigma}_{2,1}}^\ell
\lesssim& \| (a+aI(a))\Delta\vartheta^\ell\|_{\dot B^{\sigma}_{2,1}}^\ell\nonumber\\
\lesssim& \|a^\ell\Delta\vartheta^\ell\|_{\dot B^{\sigma}_{2,1}}^\ell+\|a^h\Delta\vartheta^\ell\|_{\dot B^{\sigma}_{2,1}}^\ell+\| aI(a)\Delta\vartheta^\ell\|_{\dot B^{\sigma}_{2,1}}^\ell\nonumber\\
\lesssim& \|\vartheta^\ell\|_{\dot B^{\frac 52}_{2,1}}\|a^\ell\|_{\dot B^{\sigma}_{2,1}}+
\|a^h\|_{\dot B^{\frac 3p}_{p,1}}\|\vartheta\|_{\dot B^{\sigma}_{2,1}}^\ell
+\|a\|_{\dot B^{\frac 3p}_{p,1}}^2\|\vartheta\|_{\dot B^{\sigma}_{2,1}}^\ell\nonumber\\
\lesssim&\y(t)(1+\z(t))\|(a,\vartheta)\|^{\ell}_{\dot{B}^{\sigma}_{2,1}}.
\label{jiang18}
\end{align}

The analysis on the part  $I(a){\Delta\ta^h}$ involving the high frequency $\ta^h$ is more complicated. We need a new product law different  to \eqref{ci1} and \eqref{ci2}.
To do so, we  use the  Bony decomposition
\begin{align}\label{new33}
I(a){\Delta\ta^h}= \dot{T}_{{\Delta\ta^h}}I(a)+\dot{R}({\Delta\ta^h},I(a))+\dot{T}_{I(a)}{\Delta\ta^h}
\end{align}
and  Lemmas  \ref{bernstein} and \ref{fangji}   to estimate  the first two items of right-hand side of the previous equation:
\begin{align}
\|\dot{T}_{{\Delta\ta^h}}I(a)+\dot{R}({\Delta\ta^h},I(a))\|_{\dot B^{\sigma}_{2,1}}^\ell \lesssim&\|\dot{T}_{{\Delta\ta^h}}I(a)+\dot{R}({\Delta\ta^h},I(a))\|_{\dot B^{\frac {3}{2}-\frac6p}_{2,\infty}}\nonumber\\
\lesssim&\|\dot{T}_{{\Delta\ta^h}}I(a)+\dot{R}({\Delta\ta^h},I(a))\|_{\dot B^{0}_{ p/2,\infty}}\nonumber\\
\lesssim& \|I(a)\|_{\dot B^{2-\frac 3p}_{p,1}}\|{\Delta\ta^h}\|_{\dot B^{\frac 3p-2}_{p,1}}\nonumber\\
\lesssim& \|a\|_{\dot B^{2-\frac 3p}_{p,1}}\|{\Delta\ta^h}\|_{\dot B^{\frac 3p-2}_{p,1}}\nonumber\\
\lesssim&( \|a^\ell\|_{\dot B^{\frac 12}_{2,1}}+\|a^h\|_{\dot B^{\frac 3p}_{p,1}})\|{\Delta\ta^h}\|_{\dot B^{\frac3p-2}_{p,1}}.
\end{align}

For the last term on the right-hand side of \eqref{new33}, we notice that
\begin{align*}
&\dot{\Delta}_{k} {\Delta\ta^h}\equiv0 \,\mbox{ for } k<j_0-1, \,\,\,\\
&\dot\Delta_{j}(\dot S_{k-1}I(a)\,\dot{\Delta}_{k} {\Delta\ta^h})\equiv0 \, \mbox{ for } |k-j|>4.
\end{align*}
Hence, it  follows from  Lemma \ref{bernstein}  that
\begin{align*}
\|\dot{T}_{I(a)}{\Delta\ta^h}\|^\ell_{\dot{B}^{\frac {3}{2}-\frac6p}_{2,1}}
&\lesssim \sum_{j \le j_0}2^{j(\frac32-\frac6p)}\big\|\dot{\Delta}_j\big(\sum_{k\ge j_0-1}\dot{S}_{k-1}I(a)\dot{\Delta}_{k}{\Delta\ta^h}\big)\big\|_{L^2}\\
&\lesssim \sum_{j \le j_0}2^{j(\frac32-\frac6p)}\sum_{j_0-1\le k\le j_0+4}\big\|\dot{\Delta}_j\big(\dot{S}_{k-1}I(a)\dot{\Delta}_{k}{\Delta\ta^h}\big)\big\|_{L^2}
\\&\lesssim\sum_{j \le j_0}2^{j(\frac32-\frac6p)}\sum_{j_0-1\le k\le j_0+4}(\sum_{k'\leq k-2}\|\dot{\Delta}_{k'} I(a)\|_{L^{\frac{2p}{p-2}}})\|\dot{\Delta}_k{\Delta\ta^h}\|_{L^p}
\\&\lesssim\sum_{j \le j_0}2^{j(\frac32-\frac6p)}\sum_{j_0-1\le k\le j_0+4}(\sum_{k'\leq k-2}2^{\frac{3}{p}k'}\|\dot{\Delta}_\ell I(a)\|_{L^{2}})\|\dot{\Delta}_k{\Delta\ta^h}\|_{L^p}
\\&\lesssim \|I(a)\|^\ell_{\dot{B}^{\frac 12}_{2,1}}\sum_{j \le j_0}2^{j(1-\frac{3}{p})}\sum_{j_0-1\le k\le j_0+4}\|\dot{\Delta}_{k}{\Delta\ta^h}\|_{L^p}
\\&
\lesssim \|a\|^\ell_{\dot{B}^{\frac 12}_{2,1}}\|{\Delta\ta^h}\|_{\dot{B}^{\frac 3p-2}_{p,1}}.
\end{align*}
Hence, we have the product law
\begin{align}\label{jiang19}
\|I(a)\Delta\vartheta^h\|_{\dot B^{\sigma}_{2,1}}^\ell
\lesssim \bigl(\|a\|^\ell_{\dot B^{\frac {1}{2}}_{2,1}}+\|a\|^h_{\dot B^{\frac {3}{p}}_{p,1}}\bigr)\|\vartheta\|^h_{\dot B^{\frac{3}{p}}_{p,1}}\lesssim\y(t)\z(t).
\end{align}
Therefore, by \eqref{jiang18} and \eqref{jiang19}, we have
\begin{align}\label{jiang20}
\| I(a)\Delta\vartheta\|_{\dot B^{\sigma}_{2,1}}^\ell
\lesssim&
\y(t)\z(t)+\y(t)(1+\z(t))\|\vartheta\|^{\ell}_{\dot{B}^{\sigma}_{2,1}}.
\end{align}

Finally, we consider the last nonlinear term in  \eqref{jiang2}. As in the derivation of \eqref{jiang4}, we use  \eqref{ci1} and \eqref{ci2} to obtain that
\begin{align}\label{jiang21}
&\big\|(1+I(a))\big(2|Du|^2+(\div u)^2\big)\big\|^{\ell}_{\dot{B}^{\sigma}_{2,1}}\nonumber\\
&\quad\lesssim (1\!+\!\|I(a)\|_{\dot B^{\frac{3}{p}}_{p,1}})\big\||\nabla u|^2\big\|^{\ell}_{\dot{B}^{\sigma}_{2,1}}\nonumber\\
&\quad\lesssim (1\!+\!\|a^\ell\|_{\dot B^{\frac{1}{2}}_{2,1}}\!+\!\|a^h\|_{\dot B^{\frac{3}{p}}_{p,1}})\Big((\|u\|^\ell_{\dot B^{\frac{5}{2}}_{2,1}}\!+\!\|u\|^h_{\dot B^{\frac{3}{p}\!+\!1}_{p,1}})\| u\|^{\ell}_{\dot{B}^{\sigma}_{2,1}}
\!+\!(\| u^\ell\|_{\dot B^{\frac{3}{2}}_{2,1}}\!+\!\| u^h\|_{\dot B^{\frac{3}{p}}_{p,1}})
\| u^h\|_{\dot B^{\frac{3}{p}}_{p,1}}\Big)\nonumber\\
&\quad\lesssim (1\!+\!\|a^\ell\|_{\dot B^{\frac{1}{2}}_{2,1}}\!+\!\|a^h\|_{\dot B^{\frac{3}{p}}_{p,1}})\Big((\|u\|^\ell_{\dot B^{\frac{5}{2}}_{2,1}}\!+\!\|u\|^h_{\dot B^{\frac{3}{p}\!+\!1}_{p,1}})\| u\|^{\ell}_{\dot{B}^{\sigma}_{2,1}}\nonumber\\
&\hspace{7cm}\quad+\!
\| u^\ell\|_{\dot B^{\frac{1}{2}}_{2,1}}\| u^\ell\|_{\dot B^{\frac{5}{2}}_{2,1}}\!+\!\| u^h\|_{\dot B^{\frac{3}{p}-1}_{p,1}}\| u^h\|_{\dot B^{\frac{3}{p}\!+\!1}_{p,1}}\Big)
\nonumber\\
&\quad\lesssim
(1\!+\!\z(t))\y(t)\z(t)\!+\!\y(t)(1\!+\!\z(t))\|u\|^{\ell}_{\dot{B}^{\sigma}_{2,1}}.
\end{align}

Inserting \eqref{jiang4}, \eqref{jiang8}, \eqref{jiang9}, \eqref{jiang14}, \eqref{jiang15}, \eqref{jiang20} and \eqref{jiang21} into \eqref{jiang2}, we can get
\begin{align}\label{jiang22}
&\|(a,u,\vartheta)\|^{\ell}_{\dot{B}^{\sigma}_{2,1}}+\int_0^t\|(a,u,\vartheta)
\|^{\ell}_{\dot{B}^{\sigma+2}_{2,1}}\,ds\\
&\quad\nonumber\lesssim\|(a_0,u_0,\vartheta_0)\|^{\ell}_{\dot{B}^{\sigma}_{2,1}}\!+\!\int_0^t(1\!+\!
\z(s))\y(s)\z(s)\,ds\!+\!\int_0^t\y(s)(1\!+\!\z(s))\|(a,u,\vartheta)\|^{\ell}_{\dot{B}^{\sigma}_{2,1}}\,ds.
\end{align}
From Theorem \ref{can}, we deduce that
\begin{align}\int_0^t\y(s)(1+\z(s))\,ds+\int_0^t(1+\z(s))\y(s)\z(s)\,ds\lesssim (1+\z_0)^2\z_0 .
\end{align}
Hence,
by the Gronwall inequality and the condition $\frac32-\frac6p< \sigma\le \frac 3p$,  we obtain  from \eqref{jiang22} the desired uniform bound (\ref{new3}) or
\begin{align}\label{jiang23}
\|(a,u,\vartheta)\|^{\ell}_{\dot{B}^{\sigma}_{2,1}}\le C
\end{align}
for  $t \ge 0$, where constant $C > 0$ is dependent of $\z_0$ and $\|(a_0,u_0,\vartheta_0)\|^{\ell}_{\dot{B}^{\sigma}_{2,1}}$.

\subsection{Lyapunov-type differential inequality}
We show the following Lyapunov-type differential inequality:
\begin{align}
\frac{d}{dt}
\z(t) + C\z^{\frac{ 5-2\sigma}{ 1-2\sigma}}(t)\le 0,\,\,\, t>0\label{new6}
\end{align}
for a constant $C>0$.

Indeed,   a similar proof of   \cite[Lemmas 4.1,  4.2]{xujiang2019arxiv}
implies  the following inequality
\begin{align}\label{R-E97}
\frac{d}{dt}(\|(a,u,\vartheta)\|^{\ell}_{\dot B^{\frac 12}_{2,1}}+&\|a\|^h_{\dot B^{\frac 3p}_{p,1}}+\|u\|^h_{\dot B^{\frac 3p-1}_{p,1}}+\|\vartheta\|^h_{\dot B^{\frac 3p-2}_{p,1}})
\nonumber\\
&+C\Big(\|(a,u,\vartheta)\|^{\ell}_{\dot{B}^{\frac{5}{2}}_{2,1}}+\|(a,\vartheta)\|_{\dot B^{\frac 3p}_{p,1}}^h+\|u\|_{\dot B^{\frac 3p+1}_{p,1}}^h\Big)\le0
\end{align}
or
\begin{align}
\frac{d}{dt}
\z(t) + C\y(t)\le 0,\,\,\, t>0.\label{new66}
\end{align}
Now, for  any $ \frac32-\frac6p< \sigma<\frac 12$ and $\eta_1=\frac{4}{5-2\sigma}$,
it follows from interpolation inequality in Lemma \ref{neicha}  and (\ref{new3}) that
\begin{align*}
\|(a,u,\vartheta)\|^\ell_{\dot{B}_{2,1}^{\frac 12}}
\lesssim& \big(\|(a,u,\vartheta)\|^\ell_{\dot{B}_{2,1}^{\sigma}}\big)^{\eta_1}\big(\|(a,u,\vartheta)\|^\ell_{\dot{B}_{2,1}^{\frac 52}}\big)^{1-\eta_1}\nonumber\\
\lesssim&\big(\|(a,u,\vartheta)\|^\ell_{\dot{B}_{2,1}^{\sigma}}\big)^{\eta_1}\big(\|(a,u,\vartheta)\|^\ell_{\dot{B}_{2,1}^{\frac 52}}\big)^{1-\eta_1}\lesssim\big(\|(a,u,\vartheta)\|^\ell_{\dot{B}_{2,1}^{\frac 52}}\big)^{1-\eta_1}.
\end{align*}
or the low frequency part estimate
\begin{align}\label{A7}
\big(\|(a,u,\vartheta)
\|^\ell_{\dot{B}_{2,1}^{\frac 12}}\big)^{\frac{1}{1-\eta_1}}\lesssim \|(a,u,\vartheta)\|^\ell_{\dot{B}_{2,1}^{\frac 52}}.
\end{align}
To obtain the corresponding   high frequency part estimate, we  use   the smallness property of ${\mathcal X}$ shown in (\ref{yaya}) and the high frequency property to produce
\begin{align}\label{difang}
\big(\|u\|^h_{\dot{B}_{p,1}^{\frac 3p-1}}+\|a\|^h_{\dot{B}_{p,1}^{\frac 3p}}+
\|\ta\|^h_{\dot{B}_{p,1}^{\frac 3p-2}}\big)^{\frac1{1-\eta_1}}\lesssim \|u\|^h_{\dot{B}_{p,1}^{\frac 3p+1}}
+\|a\|^h_{\dot{B}_{p,1}^{\frac 3p}}+\|\ta\|^h_{\dot{B}_{p,1}^{\frac 3p}}.
\end{align}
The combination of  \eqref{A7} and \eqref{difang} together with the definition of the symbols $\y$ and $\z$ in Theorem \ref{can} shows that    we have the lower bound estimate for the integrand in the previous equation
\begin{align}\label{new5}\z^{\frac{ 5-2\sigma}{ 1-2\sigma}}(t)=\z^\frac1{1-\eta_1}(t)\lesssim \y(t),
\end{align}
which  yields  the desired  Lyapunov-type differential inequality \eqref{new6}.

\subsection{Decay estimate}
Dividing (\ref{new6}) by $\z^{\frac{ 5-2\sigma}{ 1-2\sigma}}(t)$ and integrating the resultant inequality, we have
$$\z^{1-\frac{ 5-2\sigma}{ 1-2\sigma}}(t)\ge \z^{1-\frac{ 5-2\sigma}{ 1-2\sigma}}(0)+ (\frac{ 5-2\sigma}{ 1-2\sigma}-1)C t,$$
that is,
\begin{align*}
\z(t)\le  (\z_0^{-\frac{4}{1-2\sigma}}+\frac{4}{1-2\sigma}t)^{-\frac{1-2\sigma}{4}}
\lesssim (1+t)^{-\frac{1-2\sigma}{4}},
\end{align*}
which implies
\begin{align}\label{A10}
\|(a,u,\vartheta)\|_{\dot{B}_{p,1}^{\frac 3p-1}}\lesssim \|(a,u,\vartheta)\|^{\ell}_{\dot B^{\frac 12}_{2,1}}+\|a\|^h_{\dot B^{\frac 3p}_{p,1}}+\|u\|^h_{\dot B^{\frac 3p-1}_{p,1}}+\|\vartheta\|^h_{\dot B^{\frac 3p-2}_{p,1}}\lesssim (1+t)^{-\frac{1-2\sigma}{4}}.
\end{align}

For $\frac 3p-\frac 32+\sigma\le\alpha\le\frac 3p-1,$ by the interpolation inequality we have
\begin{align*}
\|(a,u,\vartheta)\|^\ell_{\dot{B}_{p,1}^{\alpha}}
\lesssim &\|(a,u,\vartheta)\|^\ell_{\dot{B}_{2,1}^{\alpha+\frac 32-\frac 3p}}\\
\lesssim &\big(\|(a,u,\vartheta)\|^\ell_{\dot{B}_{2,1}^{\sigma}}\big)^{\eta_{2}} \big(\|(a,u,\vartheta)\|^\ell_{\dot{B}_{2,1}^{\frac 12}}\big)^{1-\eta_{2}},\quad \eta_{2}=\frac{\frac 3p -1-\alpha}{\frac12-\sigma}\in [0,1],
\end{align*}
which combines \eqref{jiang23} with \eqref{A10} gives
\begin{align}\label{A10234}
\|(a,u,\vartheta)\|^\ell_{\dot{B}_{p,1}^{\alpha}}
\le C(1+t)^{-\frac{3}{2}(\frac 12-\frac 1p)-\frac{\alpha-\sigma}{2}}.
\end{align}
In  the light of
$\frac 3p-\frac 32+\sigma\le\alpha\le\frac 3p-1,$
 we see that
$$\|(a,u)\|^h_{\dot{B}_{p,1}^{\alpha}}\le C(\|a\|^h_{\dot{B}_{p,1}^{\frac 3p}}+\|u\|^h_{\dot{B}_{p,1}^{\frac 3p-1}})\le C(1+t)^{-\frac{1-2\sigma}{4}},
$$
from which and \eqref{A10234} gives
\begin{align*}
\|(a,u)\|_{\dot{B}_{p,1}^{\alpha}}
\le&C(\|(a,u)\|^\ell_{\dot{B}_{p,1}^{\alpha}}+\|(a,u)\|^h_{\dot{B}_{p,1}^{\alpha}})\nonumber\\
\le& C(1+t)^{-\frac{3}{2}(\frac 12-\frac 1p)-\frac{\alpha-\sigma}{2}}+C(1+t)^{-\frac{1-2\sigma}{4}}\nonumber\\
\le& C(1+t)^{-\frac{3}{2}(\frac 12-\frac 1p)-\frac{\alpha-\sigma}{2}}.
\end{align*}
Similarly, for $\frac 3p-\frac 32+\sigma\le\beta\le\frac 3p-2,$ we can get
\begin{align*}
\|\ta\|_{\dot{B}_{p,1}^{\beta}}
\le C(\|\ta\|^\ell_{\dot{B}_{p,1}^{\beta}}+\|\ta\|^h_{\dot{B}_{p,1}^{\beta}})
\le C(1+t)^{-\frac{3}{2}(\frac 12-\frac 1p)-\frac{\beta-\sigma}{2}}.
\end{align*}

Thanks to the embedding relation
$\dot{B}^{0}_{p,1}(\R^3)\hookrightarrow L^p(\R^3)$, one infer the desired decay estimates
\begin{align*}
\|\Lambda^{\alpha} (a,u)\|_{L^p}
\le C(1+t)^{-\frac{3}{2}(\frac 12-\frac 1p)-\frac{\alpha-\sigma}{2}},\quad
\|\Lambda^{\beta} \ta\|_{L^p}
\le C(1+t)^{-\frac{3}{2}(\frac 12-\frac 1p)-\frac{\beta-\sigma}{2}}.
\end{align*}

Consequently, the proof of Theorem \ref{th2} is complete.

\bigskip
\noindent \textbf{Acknowledgement.} This work is supported by NSFC under grant numbers 11601533, 11571240.


\begin{thebibliography}{99}
\linespread{0}


\bibitem{bcd}
H.~Bahouri, J.Y. Chemin, R.~Danchin,
\newblock { {F}ourier {A}nalysis and {N}onlinear {P}artial {D}ifferential
  {E}quations}.
\newblock Grundlehren Math. Wiss. , vol. {\textbf{343}}, Springer-Verlag,
  Berlin, Heidelberg, 2011.
\newblock $\,$

\bibitem{bresch}
D. Bresch, B. Desjardins,
\newblock On the existence of global weak solutions to the Navier-Stokes equations for viscous compressible and heat conducting fluids,
\newblock {\it J. Math. Pures Appl.}, {\bf 87}, 57--90, 2007.
\newblock $\,$

\bibitem{charve} F. Charve, R. Danchin,  A global existence result for the compressible Navier-Stokes equations in the critical $L^p$ framework, {\it Arch. Ration. Mech. Anal.,}  {\bf198}, 233--271, 2010.

\bibitem{lerner}
J.Y.~Chemin, N.~Lerner,
\newblock Flot de champs de vecteurs no lipschitziens
et \'{e}quations de Navier-Stokes,
\newblock {\it J. Differential Equations}, {\bf 248}, 2130--2170, 2010.
\newblock $\,$



\bibitem{chenqionglei} Q. Chen, C. Miao,  Z. Zhang, {Global well-posedness for compressible Navier-Stokes equations with highly oscillating initial velocity},  {\it Comm. Pure Appl. Math.}, {\bf63},  1173--1224, 2010.

\bibitem{chenqionglei3} Q. Chen, C. Miao,  Z. Zhang,
Well-posedness in critical spaces for the compressible
Navier-Stokes equations with density dependent viscosities, {\it Revista Mat. Iber.,} {\bf26}, 915--
946, 2010.

\bibitem{chenqionglei2}
Q. Chen, C. Miao, Z. Zhang,  {On the ill-posedness of the compressible Navier-Stokes equations in the critical Besov spaces}, {\it Rev. Mat. Iberoam.}, {\bf31}, 1375--1402, 2015.

\bibitem{zhaixiaoping}
 Z. M.~Chen, X.~Zhai,
\newblock
{ Global large solutions and incompressible limit for the compressible Navier-Stokes equations},
\newblock {\it J. Math. Fluid Mech.,} https://doi.org/10.1007/s00021-019-0428-3, 2019.
\newblock $\,$


\bibitem{chikami} N. Chikami, R. Danchin, {On the well-posedness of the full compressible Navier-Stokes system in critical Besov spaces}, {\it J. Differential Equations},  {\bf258},  3435--3467, 2015.

\bibitem{cj} Y. Cho, B.J. Jin,  Blow-up of viscous heat-conducting compressible flows, {\it J. Math. Anal. Appl.},
{\bf320},  819--826, 2006.

\bibitem{danchin2000} R. Danchin, {Global existence in critical spaces for compressible Navier-Stokes equations}, {\it Invent. Math.}, {\bf141}, 579--614, 2000.

\bibitem{danchin2001cpde} R. Danchin, {Local theory in critical spaces for compressible viscous and heat-conductive gases}, {\it Comm. Partial Differential Equations}, {\bf26},   1183--1233, 2001.

\bibitem{danchin2001arma} R. Danchin, Global existence in critical spaces for flows of compressible viscous and heat-conductive gases,   {\it Arch. Rational Mech. Anal.}, {\bf 160}, 1--39, 2001.


\bibitem{helingbing} R. Danchin,  L. He,  The incompressible limit in $L^p$ type critical spaces,
{\it Math. Ann.}, {\bf 366}, 1365--1402,  2016.

\bibitem{xujiang2017arma} R. Danchin,  J. Xu, Optimal time-decay estimates for the
compressible Navier-Stokes equations in the critical $L^{p}$ framework, {\it Arch. Rational Mech. Anal.}, {\bf 224}, 53--90, 2017.

\bibitem{xujiang2019jmfm} R. Danchin,  J. Xu,  Optimal decay estimates in the critical $L^{p}$ framework for flows of compressible viscous and heat-conductive gases, {\it  arXiv:1612.05776}.


\bibitem{Feireisl}
E. Feireisl, {Dynamics of Viscous Compressible Fluids}. \textit{Oxford Univ. Press, Oxford,} 2004.





\bibitem{Feireisl3}
E. Feireisl, {On the motion of a viscous, compressible and heat conducting fluid}, {\it Indiana Univ. Math. J.}, {\bf53},  1705--1738, 2004.

\bibitem{Feireisl2} E. Feireisl, A. Novotn\'{y}, H. Petzeltov\'{a}, { On the global existence of globally defined
weak solutions to the Navier-Stokes equations of isentropic compressible fluids}, {\it J. Math. Fluid Mech.},  {\bf3},
358-392, 2001.


\bibitem{haspot} B. Haspot, {Well-posedness in critical spaces for the system of compressible Navier-Stokes in larger spaces}, {\it J. Differential Equations}, {\bf251}, 2262--2295, 2011.

 \bibitem{hoff} D. Hoff, H. Jenssen, {Symmetric nonbarotropic flows with large data and forces}, {\it Arch. Ration. Mech. Anal.}, {\bf173} (2004), 297--343.

 \bibitem{huangxiangdi2}
X. Huang, J. Li,  {On breakdown of solutions to the full compressible Navier-Stokes equations}, {\it Meth. Appl. Anal.}, {\bf16}, 479--490, 2009.

\bibitem{huangxiangdi4}
X. Huang, J. Li,  {Global classical and weak solutions to the three-dimensional full compressible Navier-Stokes system with vacuum and large oscillations}, {http://arxiv.org/abs/1107.4655v3 [math-ph],} 2011.

\bibitem{huangxiangdi1}
X. Huang, J. Li,  {Serrin-type blowup criterion for viscous, compressible, and heat conducting Navier-Stokes and magnetohydrodynamic flows}, {\it Comm. Math. Phys.}, {\bf324}, 147--171, 2013.


\bibitem{huangxiangdi3}
X. Huang, J. Li, Y. Wang, { Serrin-type blowup criterion for full compressible Navier-Stokes system}, {\it Arch. Ration. Mech. Anal.}, {\bf207}, ( 303--316, 2013.


  \bibitem{itaya}
N. Itaya, {On the Cauchy problem for the system of fundamental equations describing themovement of compressible viscous fluid},
{\it Kodai Math. Semin. Rep.}, {\bf23},  60--120, 1971.

\bibitem{jiangsong1998}S. Jiang, {Large-time behavior of solutions to the equations of a viscous polytropic ideal gas}, {\it Ann. Mat. Pura Appl.}, {\bf175},   253--275, 1998.

\bibitem{jiangsong1999}S. Jiang, { Large-time behavior of solutions to the equations of a one-dimensional viscous polytropic ideal gas in unbounded domains}, {\it Commun. Math. Phys.}, {\bf200}, 181--193, 1999.


\bibitem{ks}A.V. Kazhikhov, V.V. Shelukhin, {Unique global solution with respect to time of
initial-boundary value problems for one-dimensional equations of a viscous gas},  {\it J. Appl.
Math. Mech.,} {\bf41},  273--282, 1977.

\bibitem{kagei} Y.~Kagei, T.~Kobayashi,  Asymptotic behavior of solutions of the compressible
Navier-Stokes equations on the half space, {\it Arch. Rational Mech. Anal.}, {\bf{177}}, 231--330,  2005.

\bibitem{kobayashi} T.~Kobayashi,  Y. Shibata,  Decay estimates of solutions for the equations of motion
of compressible viscous and heat-conductive gases in an exterior domain of  $\mathbb{R}^3$,
{\it Com. Math. Phys.},  {\bf{200}}, 621--659,  1999.



\bibitem{ma2} A. Matsumura,  T. Nishida,
The initial value problem for the equations of motion of viscous and
heat-conductive gases, {\it  J.  Math. Kyoto Univ.},
 {\bf 20},   67--104,  1980.

 \bibitem {ma3}
A. Matsumura,   T.~Nishida,  Initial boundary value problems for the equations of motion of compressible viscous and heat conductive fluids,
{\it Commun. Math. Phys.}, {\bf{89}}, 445--464,  1983.












\bibitem{nash} J. Nash, {Le probl\`eme de Cauchy pour les \'equations diff\'erentielles d'un fluide g\'en\'eral},
{\it Bulletin de la Soc. Math. de France,} {\bf90}, 487--497, 1962.



\bibitem{ponce} G.~Ponce,  Global existence of small solution to a class of nonlinear evolution equations,
{\it Nonlinear Anal. TMA.}, {\bf{9}}, 339--418,  1985.


\bibitem{tani}
A. Tani, { On the first initial-boundary value problem of compressible viscous fluid motion}, {\it Publ. Res. Inst. Math. Sci. Kyoto Univ.}, {\bf13},  193--253, 1977.


\bibitem{wenhuanyao}
H. Wen,  C. Zhu,
{Global symmetric classical solutions of the full compressible Navier-Stokes equations with vacuum and large initial data,}
 {\it J. Math. Pures Appl.},  {\bf102},  498--545, 2014.



\bibitem{xinzhouping} Z. Xin,
{Blow up of smooth solutions to the compressible Navier-Stokes equation with compact
density}, {\it Comm. Pure Appl. Math.}, {\bf51},   229--240, 1998.

\bibitem{xujiang2019arxiv}
Z. Xin, J. Xu,
{Optimal decay for the compressible Navier-Stokes equations without additional smallness assumptions},
\newblock {\it  ArXiv:1812.11714v2}.
\newblock $\,$

\bibitem{xujiang2019cmp}
J. Xu,
{A low-frequency assumption for optimal time-decay estimates to the compressible Navier-
Stokes equations},
\newblock {\it Comm. Math. Phys.,  https://doi.org 10.1007/s00220-019-03415-6, 2019}.
\newblock $\,$


\end{thebibliography}
\end{document}